\documentclass{article}
\usepackage{amsmath,amssymb,latexsym,amsthm}
\font\goth=eusm10

\newcommand\M{\mathcal M}
\newcommand\E{\mathcal E}
\newcommand\Ii{\hbox{\goth I}}

\newcommand\Oc{\hbox{\goth O}}

\newcommand\CC{\mathbf C}
\newcommand\RR{\mathbf R}

\newcommand\Pd{\mathbf{P}^2}
\newcommand\Pj{\mathbf{P}}

\newcommand\Pc{\mathbf{P}^5}
\newcommand\Pt{\mathbf{P}^3}
\newcommand\ZZ{\mathbf{Z}}
\newcommand\QQ{\mathbf{Q}}

\newtheorem{proposition}{Proposition}[section]

\newtheorem*{definition}{Definition}

\newtheorem{theorem}{Theorem}[section]

\newtheorem*{problem}{Problems}

\newtheorem{example}[theorem]{Example}

\theoremstyle{plain}
\theoremstyle{definition}
\theoremstyle{remark}
\begin{document}
\title{A remark on K3s of Todorov type $(0,9)$ and $(0,10)$}
\author{C.Madonna}
\date{}
\maketitle
\begin{abstract}
Frequentely it happens that isogenous (in the sense of Mukai)
K3 surfaces are partners of each other and sometimes
they are even isomorphic.
This is due, in some cases, to the (too high, e.g. 
bigger then or equal to 
$12$)) rank of the Picard lattice
as showed by Mukai in \cite{Muk3}.
In other cases this is due to the structure of the Picard
lattice and not only on its rank. This is the case, for example, 
of K3s with Picard lattice containing a latice of Todorov type $(0,9)$ or $(0,10)$.
\end{abstract}
\section{Introduction and main statements}

In recent times it seems that the following lines, written
on page 394 (see lines 11--12) in \cite{Muk3} 
by Mukai during the 80s, received more and more attention:\par
\smallskip
{\it "...there is a pair of K3 surfaces $S$ and $M$ such that
$T(S) \cong T(M)$, but $N(S) \not\cong N(M)$, as lattices."} \par
\smallskip
To start with let us recall the following:
\begin{definition}
Let $S$ be a K3 surface. A K3 surface $M$ is called
a Fourier--Mukai partner of $S$ if there exists
an Hodge isometry $(T(S),\CC \omega_S) \to (T(M),\CC \omega_M)$
(i.e. an isomorphism of lattices preserving the Hodge decomposition),
where $\omega_S,\omega_M$ are the nowhere vanishing holomorphic two forms
on $S$ and $M$, respectively.
\end{definition}

Here, by a K3 surface we mean a smooth projective surface $S$ over $\CC$
with $\Oc_S(K_S) \cong \Oc_S$ and $h^1(\Oc_S)=0$, and 
we consider the Hodge decomposition of $H^*(S,\CC)=H^*(S,\ZZ) \otimes \CC$
given by
\[
H^{*(i,j)}(S,\CC)=H^{(i,j)}(S,\CC) \ {\rm for} \
(i,j) \in \{ (2,0),(0,2) \}
\]
\[
H^{*(1,1)}(S,\CC)=H^0(S,\CC) \oplus H^{(1,1)}(S,\CC) \oplus H^4(S,\CC).
\]
It is well known
that there exists pairs of $K3$ surfaces
$S$ and $M$ with Hodge isometric trascendental lattices but non isomorphic
Picard groups. 
In particular these pairs are given by non isomorphic
K3s. By the above, it is clear that the Hodge isometry class of 
the trascendental lattice $T(S)$
does not recover neither the K3 surface $S$ nor the Picard lattice,
but only the genus of the Picard lattice, and hence
the interest for K3s with Picard
lattice containing only one class in its genus, naturally arises.
For, in this case, one can still hope to recover the K3 itself.
Nowadays, one believes in the existence of such pairs, 
from a lattice theoretical point of view as suggested in
Nikulin's paper \cite{nikulin} 
(see also page 394 in \cite{Muk3}),
from a categorial point of view as in Bridgeland and Maciocia's paper
(see the proof of Proposition 5.3 in
\cite{bridgeland-maciocia}) or  Orlov's paper (see
\cite{orlov-abelian}, the example on the last page), 
but also from a geometric point of view
for which we refer to Morrison's paper (see page 10 in \cite{morrison}).

Recently, Oguiso in \cite{oguiso} showed that for any integer
$n$ one can give $n$ pairs of K3s as above.
Moreover, given a K3 surface of Picard number $\rho=1$,
he computed the number of non isomorphic K3 surfaces 
(of Picard number equal to 1)
with Hodge isometric trascendental lattices, i.e.
the number of non isomorphic Fourier--Mukai partners of the given K3.
In other words such a pair, say $S$ and $M$,
is characterized by the existence
of a Hodge isometry between the relative trascendental lattices which
does not extend to an Hodge isometry between the K3 lattices
$H^2(S,\ZZ)$ and $H^2(M,\ZZ)$.
On the other hand, as showed by Mukai in \cite{Muk3}, 
one can always extend such an isometry
between the relative Mukai lattices $\tilde{H}(S,\ZZ)$ and 
$\tilde{H}(M,\ZZ)$.
Let us recall the
Mukai lattice is given by the group $H^*(S,\ZZ)=H^0(S,\ZZ) \oplus
H^2(S,\ZZ) \oplus H^4(S,\ZZ)$ endowed with
the Mukai product, defined for vectors $v=(v_0,v_1,v_2)$ and
$w=(w_0,w_1,w_2)$ by setting $(v,w)=v_1 \cdot w_1-v_0 \cdot w_2-v_2 \cdot w_0$,
i.e. a Hodge structure of weight $2$ on the cohomology ring of $S$, defined by mean
of the cup product $\cdot$.
The obstruction to extend the isometry between $T(S)$ and $T(M)$
to the relative K3 lattices is also related
to the existence of non isomorphic primitive embeddings of 
$T \hookrightarrow \Lambda$, where $\Lambda$ is the abstract K3
lattice.
On the other hand, 
this is not the case for K3s of Picard number $\rho \geq 12$ for which
the primitive embedding $T \hookrightarrow \Lambda$ is always unique.
From a lattice point of view, 
one could look for other conditions to ensure such
type of unicity result, e.g. under some conditions on the lattice
or on the orthogonal to its image under such an embedding.
In this direction we quote
a result of Morrison (see \cite{morrison-todorov} Theorem A1)
which is the key point of our proof of Theorem \ref{thm:todorov}
for K3s of Todorov type $(0,10)$ and Picard number equal to $11$,
and a computation by Miranda and Morrison of the
number of equivalence classes of such primitive embeddings 
$\varphi:T \hookrightarrow \Lambda$ (see 
\cite{miranda-morrison}).
Later this computation was generalized in \cite{HLYO2} to get a formula
to compute the number of Fourier--Mukai partners of any K3 surface. 
Note, 
the relevance of both the knowledge of the genus of a lattice $N$ and 
the 
map between the isometries groups of $N$ and of its discriminant 
$A_N:=N^{\vee}/N$,
for the case of the Picard lattice 
$NS(S)=N$ of a K3 surface $S$, is the main ingredient if
one wants to compute the number of Fourier--Mukai
partners of $S$ by mean of this formula. Unfourtunately, as far as we know,
there is not a general algorithm
to compute the genus of a lattice and hence 
the computation by mean of such a formula can be not so trivial.
On the other hand if one can find some geometric argument
for such type of computation, one can use, in many cases,
the above formula, to get the genus of the lattice.
For example if $S$ is a K3 surface for which the map
$O(N') \to O(A_{N'})$ is surjective for all $N'$ in the genus $g(N)$ of
$N=NS(S)$, then the number of classes contained in $g(N)$ is given by the number of
Fourier--Mukai partners of $S$.
We refer to \cite{cassels} or \cite{nikulin}
for conditions 
to ensure both the existence of only one class in the genus
of a lattice $N$ and the surjectivity of the map $O(N) \to O(A_{N})$.
In this direction some computations are given  
in \cite{HLYO1} \cite{HLYO2} \cite{madonna-retta} \cite{oguiso}, for example.

In this note we consider two special cases on the borded, with respect
to the Picard number, i.e. K3 surfaces of Picard number 
equal to $11$ and $10$, and we focus on the following:

\begin{problem} \label{problem:pro} \ \
\begin{enumerate}
\item[a.] Let $S$ and $M$ be two polarized K3 surfaces with
Hodge isometric trascendental
lattices, $T(S) \cong T(M)$. 
Find conditions to ensure that $S \cong M$.
\item[b.] As above replacing the
Hodge isometry with an isomorphism of Hodge structures
preserving the Hodge decomposition:
\[
T(S) \otimes \QQ \to T(M) \otimes \QQ.
\]
\end{enumerate}
\end{problem}

Problem a. is related with K3s for which the Fourier--Mukai 
partner is unique
and we refer the reader to \cite{bridgeland-maciocia}
\cite{oguiso} \cite{madonna-retta}, for example.
Problem b. concerns on
conditions such that an isomorphism of Hodge structure
preserving the Hodge decomposition,
is induced by an isomorphism, a special case of a more general
phenomena regarding isomorphism induced by compositions of prime
isogenies and elements representing zero in the Picard lattice
of a K3 (see \cite{nikulin-zero}).
Let us remark there exist "many" examples of K3 surfaces
for which any isomorphism of Hodge sructure
is induced by an isomorphism between the K3s.
As remarked this is the case of isogenous (in the sense of Mukai)
K3 surfaces with Picard number bigger then or equal to 12.
\par
\smallskip
Clearly, there are also many examples of isogenous K3 surfaces
$(S,M)$ which are not isomorphic, and moreover with not Hodge isometric
trascendental lattices, but always isomorphic Hodge structure
(i.e. over $\QQ$).
In this case we are dealing for example, as showed by Caldararu in
\cite{caldararu}, with non fine moduli spaces $M$ of sheaves over $S$.
The non fineness property is related, in this case, to the existence
of an element 
in the Brauer group of the moduli space
$\alpha \in Br(M)$ (see \cite{caldararu})
or, equivalently, to the existence of a non trivial element
in the cokernel of the map $T(M) \to T(S)$ induced by the map
on the Mukai lattices constructed by means of a quasi universal sheaf
on the product $S \times M$ (see \cite{Muk3}).
In this case 
the example we have in mind is the one given by Mukai in \cite{Muk3},
i.e. $S$ is the base locus of a general net of quadrics in $\Pj^5$,
and $M$ is a sextic double plane of Picard number equal to 1 -- the moduli
space of $H$--stable rank--two sheaves over $S$ with Chern classes 
$c_1=H$ and $c_2=4$, being $H$ the class of an hyperplane section 
of $S$ (see Example 2.8 in \cite{Muk3}).
\par
\smallskip
In this note we consider a special case
of the border case $\rho(S)=11$, which is in some way our
working example.
For, let us consider the K3 cover $S \to T$ of an Enriques surface $T$.
It was showed by Verra (see \cite{Verra2}) the existence
on $S$ of a polarization $H$ of degree $8$, i.e. $H^2=8$, 
such that $S$ is the base locus of a net of quadrics in $\Pj^5$
and the 
involution induced by the natural involution of $\Pj^5$ is such that $S/ i \cong T$.
Moreover it was showed by Cossec (see \cite{cossec}), 
when $T$ is of special type, say a Raye Congruence, then there
exists $h$ on $S$ of degree $2$ such that $S \to \Pj^2$
is a double cover ramified along the union of two generically smooth cubics
intersecting transversally endowed with a totally tangent
conic out of the intersection points of the two cubics.
We focus our attention on these two cases.
First, we relate the two classes $h$ and $H$,
using the same procedure used by Mukai for the general case.
This is done considering the K3 cover $S$ as the base locus
of a special net of quadrics in $\Pj^5$ as in \cite{Verra2}.

\begin{theorem} \label{thm:reye}
Let $S$ be a K3 surface which is universal 
covering of an Enriques surface $T$.
Then, there exists $H$ on $S$ of degree $8$ such that
the moduli space $\M_S(2,H,2)$ is a double plane ramified
along the union of two generically smooth cubics, which are
fixed by the covering involution.
When $T$ is of type Reye congruence, then there exists an Enriques--Cayley
class $h$ of degree $2$ over $\M_S(2,H,2)$.
In particular $S$ and $\M_S$ are isomorphic.
\end{theorem}

The second statement of the previous theorem is a consequence,
after having showed that $\M_S$ is a partner of $S$, by the 
of the following more general result
on (the minimal desingularization of) K3 surfaces of Todorov type, which is
our second main remark:

\begin{theorem} \label{thm:todorov}
The number of Fourier--Mukai partners of a K3 surface of Todorov type 
is equal to $1$.
\end{theorem}

The proof of previous theorem is essentially contained in the proof
of Theorem 6.1 in \cite{morrison-saito} if one follows
the same strategy 
as the one used by Mukai in \cite{Muk3} for the
case of K3 surfaces of Picard number bigger then or equal to $12$.
Note also, our result is not so surprise since the Picard number
of a K3 of Todorov type is an integer between 10 and 17, and from 12 on
the result was yet known. So our contribution regards only the cases for which
$\rho=10$ and $\rho=11$.
We advise the reader our notion of K3 surface of Todorov type
does not include the ones of type $(5,16)$.
This case will be considered in Example 3.1 as an application
of the formula derived in \cite{HLYO2}.

Finally, few words on Theorem \ref{thm:reye}.
In general $\M_S$ is not a partner of $S$, unless it is isomorphic to
$S$, but one can still
consider it as an $\alpha$--twisted partner in the sense of Caldararu
(see \cite{caldararu}).

{\bf Aknowledgement.} I have to thank Professor A.Verra for his 
interest in the subject o this paper.

\section{A proof of Theorem \ref{thm:todorov}}

In this section we give a proof of Theorem \ref{thm:todorov}.
Let us briefly recall some generalities on Todorov lattices
and K3s of Todorov type.
We refer the reader to \cite{morrison-todorov} or \cite{morrison-saito}
for more details.

\begin{definition}
Let $\Lambda$ be the K3 lattice and let $E_k$ be the lattice
generated by the elements $e_1,\cdots,e_k$ with $e_i \cdot e_j=-2\delta_{ij}$,
$\delta_{ij}=1$ if $i=j$ and $\delta_{ij}=0$ otherwise.
Let $j:E_k \hookrightarrow \Lambda$ be an embedding.
The "double point lattice $L_{\alpha,k}$" is the lattice
defined as the saturation of $E_k \otimes \QQ \cap j^{-1}(\Lambda)$.
We have $[L_{\alpha,k}:E_k]=2^{\alpha}$ 
and $\alpha$ is called the ``2--index of $L_{\alpha,k}$''.
The "Todorov lattice $M_{\alpha,k}$" is the lattice
generated by $L_{\alpha,k}$ and the elements $\lambda, \mu$ such that:
\begin{enumerate}
\item[(i)] $\alpha \leq 4$, $k \geq 9$, and $2^{4-\alpha}(2^{\alpha}-1) 
\leq k \leq \alpha +11$
with $(\alpha,k) \ne (1,9)$;
\item[(ii)] $\lambda^2=2k-16$;
\item[(iii)] $\lambda \cdot L_{\alpha,k}=0$;
\item[(iv)] $\mu=\frac12 (\lambda+\sum_i e_i)$.
\end{enumerate}
\end{definition}

Our main interest in this paper are K3 surfaces $S$ for which
$NS(S)$ contains a Todorov lattice, in this way we can reduce
the proof to an arithmetic condition on the Picard lattice.
Let us recall in \cite{morrison-todorov}
it is showed the existence of Todorov lattices for each of the following pairs:
\[
(\alpha,k) \in \{ (0,9),(0,10),(0,11),(1,10),(1,11),(1,12), 
\]
\[
(2,12),(2,13),(3,14),(4,15) \}.
\]
As remarked in the introduction the maximal rank of a Todorov 
lattice is $17$, while by the above list one gets $16$ as maximum.
We prefered to define Todorov lattice in such a way not including 
this last of rank $17$ and type $(5,16)$,
included in \cite{morrison-todorov},
since it corresponds to Kummer surfaces for which our
main result will be derived in Example \ref{ex:(5,16)},
as an application of the formula
for Fourier--Mukai partners derived in \cite{HLYO2}.

\begin{definition}
A "K3 surface of Todorov type $(\alpha,k)$"
is a K3 surface $S$ with rational double points
endowed with an ample line bundle
$l$ of degree $l^2=2k-16$ and a partial desingularization
$\nu:\Sigma \to S$ with $k$ ordinary double points which has 2--index
$\alpha$ and such that
if $\pi:X \to \Sigma$ and $\mu=\nu \circ \pi$ are the minimal desingularizations
then
\[
\frac12 c_1(\mu^*(l)) \otimes \Oc_X(\sum_{p \in Sing \Sigma} \pi^{-1}(p))
\in H^2(X,\ZZ).
\]
\end{definition}

For a general K3 surface $S$ 
of Todorov type $(\alpha,k)$ we have $NS(S)=M_{\alpha,k}$
and $\Sigma=S$. In this paper we will frequently refer to a K3 of Todorov type
as to its minimal desingularization.

\begin{example}
\normalfont{
An example of a K3's of Todorov type $(0,9)$ is given by
double covering of a plane branched along the union of two cubics $C_1,C_2$
intersecting transversally.
For the case $(0,10)$ we require the existence of a totally tangent 
conic to the two cubics out of the intersection
points of the two cubics.
Let us remark, as showed by Cossec (see \cite{cossec}),
the K3 cover of an Enriques surface of type Reye Congruence admits such type
of double plane representation.
As showed by Morrison and Saito
in \cite{morrison-saito}, the action of the plane Cremona transformations
on such planes provides other examples of such K3s.
}
\end{example}

Let us start recall the following:
\begin{definition}
The "signature" of a non degenerate lattice $L$ is the pair $(r_+,r_-)$
describing the number of positive and negative eingevalues 
of the induced quadratic form on $L \otimes \bf R$.
A "positive sign structure" on a non degerate lattice $L$ of signature $(r_+,r_-)$
is a choice of one of the connected components of the set of oriented
$r_+$--planes in $L \otimes \RR$
on which the form is positive definite.
The sign structure containing the oriented $\nu$ plane
is denoted by $[\nu]$.
$O(L)$ denotes the orthogonal group of $L$
and $O_-(L)$ denotes the subgroup of isometries preserving
a sing structure.
\end{definition}

\begin{definition}
Let $j:M_{\alpha,k} \hookrightarrow \Lambda$ be a primitive embedding
of $M_{\alpha,k}$ into the K3 lattice $\Lambda$ and denote by 
$T_{\alpha,k}=j(M_{\alpha,k})^{\perp}$.
Fixed a positive sign structure $[\nu_{\alpha,k}]$ on 
$T_{\alpha,k}$
the "period space" ${\bf D}_{\alpha,k}$
is defined by setting
\[
{\bf D}_{\alpha,k}=
\{ \omega \in \Pj(T_{\alpha,k} \otimes \CC): \omega \cdot \omega
=0, \omega \cdot \overline{\omega} >0, Re(\omega) \wedge Im(\omega) \in [\nu_{\alpha,k}] \}.
\]
\end{definition}

We also need the following (see \cite{Muk3}):
\begin{definition}
Two K3 surfaces $S$ and $M$ are called ``isogenous 
in the sense of Mukai'' if there exists
an algebriac cycle $Z \in H^4(S \times M,\QQ)$
inducing an isometry between $H^2(S,\QQ)$ and $H^2(M,\QQ)$
by setting $t \to \pi_{M,*}(Z \cdot \pi^*_S(t))$.
In particular this is equivalent to the existence
of an algebraic cycle $Z$ on $S \times M$ such that
the component $Z_t$ of $Z$ in $T(S) \otimes T(M) \otimes \QQ$
induces a rational cohomological isogeny between $S$ and $M$,
i.e. an isometry of forms $T(S) \otimes \QQ \to T(M) \otimes \QQ$.
\end{definition}

Let us remark isogenous K3 surfaces in the above sense
are not always partners each other.

Let $S$ be a K3 surface of Todorov type and $M$ be a partner of $S$.
Let us denote by
$\varphi:T(S) \to T(M)$ the induced Hodge isometry
between the relative trascendental lattices.
To show Theorem \ref{thm:todorov} our main tools are the following results
(see Theorem 3.3 and
6.1 in \cite{morrison-saito}, respectively):

\begin{proposition} \label{prop:nikulin-type}
Let $M_{\alpha,k}$ be a Todorov lattice and $f:M_{\alpha,k} \hookrightarrow \Lambda$
be a primitive embedding into the K3 lattice. Let us denote by
$T_{\alpha,k}=f(M_{\alpha,k})^{\perp}$ the orthogonal complement of the image of $f$.
If $\psi_1,\psi_2:T_{\alpha,k} \hookrightarrow \Lambda$ are two
primitive embeddings then there is some $\gamma \in O_-(\Lambda)$ such that
$\gamma \circ \psi_1=\psi_2$.
\end{proposition}

\begin{proposition} \label{prop:torelli-type}
Two K3 surfaces of Todorov type, say $S$ and $M$,
giving the same point in the period space of K3s of type $M_{\alpha,k}$
are birational isomorphic.
Moreover there is an isomorphism between the minimal desingularizations
$\tilde{S}$ and $\tilde{M}$ of $S$ and $M$, respectively.
\end{proposition}
\par
\smallskip
{\it Proof of Theorem \ref{thm:todorov}.}
Let $M$ be a partner of $S$ and
consider the two primitive 
embeddings $T(S) \hookrightarrow \Lambda$ and $T(M) \hookrightarrow
\Lambda$. Since as lattices $H^2(S,\ZZ) \cong H^2(M,\ZZ) \cong \Lambda$
by Proposition \ref{prop:nikulin-type} there exist $\gamma \in O_-(\Lambda)$
such that $\gamma_{|T(S)}=\varphi$. The conclusion follows now
by Proposition \ref{prop:torelli-type}.

\section{A proof of Theorem \ref{thm:reye}}

Let $S$ be a 2--dimensional complete intersection of three quadrics $Q_1,Q_2,Q_3$ 
in $\Pj^5$ and $N$ be the net spanned by the quadrics $Q_i$.
The discriminant of such a net is a plane sextic, generically smooth,
representing the locus of singular quadrics in the net.
There are some special nets which are our main interest in this section
being related to K3 surfaces universal covering of Enriques surfaces, which we recall
in the following:
\begin{definition}
A net of quadrics $N$ in $\Pj^5$ is called an "Enriques net
of quadrics" if the base locus $Bs(N)$ is a K3 surface
universal covering of an Enriques surface $T$.
\end{definition}

Generically a net of such type contains nine rank 4 quadrics and the 
discriminant curve
splits as the union of two smooth plane cubics intersecting transversally
(see \cite{Verra2}).
The double covering of the net ramified along the union of two
cubics as above, as remarked, is a K3 of Todorov type $(0,9)$.
Note the cubics are fixed by the covering involution.

Let us consider the base locus $S$ of a general net of quadrics $N$ in $\Pj^5$ and
denote by $\M_S$ the moduli space of $H$--stable rank--two sheaves
with Chern classes $c_1=H$ and $c_2=4$, i.e. with Mukai vector
$v=(2,H,2)$. 
Mukai showed (see \cite{Muk3}) the following:

\begin{proposition} 
The moduli space $\M_S$ is a K3 double plane ramified along a smooth
plane sextic -- the double cover of the net $N$ ramified along
the discriminant locus of the net.
\end{proposition}

\begin{proof} 
Let $S$ be the base locus of a general net $N$ of quadrics in $\Pj^5$.
Let us recall that
if $\E=[\E] \in \M_S$ then we have an exact sequence
\[
0 \to \Oc_S \to \E \to \Ii_Z(H) \to 0
\]
and
$h^1(\E)=h^1(\Ii_Z(H))-h^2\Oc_S=h^1(\Ii_Z(H))-1$.
By the standard exact sequence
\[
0 \to \Ii_Z(H) \to \Oc_S(H) \to \Oc_Z(H) \to 0
\]
since $h^0(\Oc_S(H))=6$ it follows, for general $Z$,
$h^1(\Ii_Z(H))=1$, $h^1(\E)=0$ and hence $h^0(\E)=4$.
Moreover since
$S$ is complete intersection of three quadrics in
$\Pc$, and since we have assumed $S$ "general", then
$S$ is the base locus of a net of quadrics
$N$ containing quadrics of rank bigger then or equal
to $5$.
For such a net it is canonically defined a double covering map
$\pi:\tilde{N} \to N$
which is generically $2:1$; for, on a smooth quadric $Q$ the
preimage
$\pi^{-1}(Q)$ consists of the two rulings of $Q$. Clearly, this
map
ramifies along the set of quadrics for which there exists only one
family
of maximal linear subspaces, and hence it ramifies along the set
of singular quadrics $N_o$. We will show that
$\tilde{N}=\M_S$.
It is known that $\dim \M_S=2$.
For, since $h^0(\E)=4$ for a general $\E \in \M_S$ then the section
$Z_t=\{ s_t=0: t \in \Gamma \}$, $\Gamma \cong
\Pt$, moves in a 3-dimensional family of sections
which in turn correspond to one of the two 3-dimensional family of planes
$\Delta=\{ \Pd_t: t \in \Gamma \}$ contained in
a smooth quadric
$Gr(2,4) \subset \Pc$ simply setting
$Z_t=\Pd_t \cap S$.
For, bundles coming from smooth quadrics arise in pairs, i.e. to a quadric
$Q \in N$
corresponds two rank--two bundles fitting in the following
exact sequence
\[
0 \to \overline{\tau} \to \Oc_{Q}^4 \to \tau \to 0,
\]
where $\overline{\tau}$ and $\tau$ are the universal subbundle and
quotient bundle on $Q=Gr(2,4) \subset \Pj^5$, and define by
restriction two bundles $\overline{\E}=(\overline{\tau}_{\mid S})^{\vee}$ 
and $\E=\tau_{\mid S}$ on $S$. Clearly
$c_1(\E)=c_1(\overline{\E})$ and
$c_2(\overline{\E})=c_1^2-c_2(\E)=4$. Hence it is defined a 2:1
map $\M_S \to N$ and this map ramifies along
the set of quadrics containing $S$ for which there exists
only one family of maximal linear subspaces and hence along
singular quadrics containing $S$, and $\tilde{N}=\M_S$.
\par\noindent
Since the set of
singular quadrics in $\Pc$ is a hypersurface of degree $6$ in the
space of all quadrics in $\Pc$, $\vert \Oc_{\Pc}(2)
\vert=\Pj^{20}$, and since the quadrics of the net are not all
singular, then
the set of singular quadrics containing $S$ is a
hypersurface $N_o$ of degree $6$ in the plane $N=\Pj^2$:
it is a plane sextic
curve $N_o$.
Let us show that it is smooth. First, since we
assumed that all quadrics containing $S$ have rank bigger than or equal
to $5$, then the set of singular quadrics containing $S$ is
given by the set of rank--5 quadrics. Note that, since $S$ is
smooth, $Q$ is a smooth point for $N_o$ if and only if ${\rm
Sing}(Q)$ is a point. If $Q \in
N_o$ then $Q$ is a quadric in $\Pc$ with ${\rm rk}(Q)=k \leq 5$
hence ${\rm Sing}(Q)=\Pj^{5-k}$ is a point if and only if $k=5$.
Since in our hypotheses, every quadric containing $S$ has
rank bigger than or equal to $5$ every point $Q \in N_o$ is a
smooth point and the plane sextic $N_o \subset N$ is
smooth.\par\noindent
In particular it follows that $\M_S$
is again a K3, double covering of $\Pd=N$, $\phi:\M_S
\to \Pd$, branched along a smooth plane sextic curve $N_o \subset
\Pd$.
Note that while $S$ is a K3 surface of
genus $5$ and degree $8$, $\M_S$ is a K3 surface of genus
$2$ and degree $2$ on which we have the
natural involution $i$ which interchanges the sheets
of the cover, since it interchanges the rulings of the quadric on which acts;
the locus of fixed points of such involution, i.e
the subset of points $x \in \M_S$ such that
$i(x)=x$, is clearly the preimage $\phi^{-1}(N_o)$, since
only for this last there
exist only one family of maximal linear subspaces.
For, if $x \in N \setminus N_o$, then $x$ represent a rank--6
quadric $Q$ containing $S$ and the two families $\Delta$ and
$\overline{\Delta}$ of maximal linear subspaces contained in $Q$
represent the global sections $Z_t=\Pd_t \cap S$ and
$\overline{Z}_t=\overline{\Pd}_t \cap S$, of two rank--two bundles
$\E,\overline{\E}$ on $S$. Since the elements of the two family as
above are interchanged by the involution $i_S$ also the section
$Z_t$ and $\overline{Z}_t$ are interchanged and hence the bundles
represented by the above quadrics are also interchanged, and any
one of them is the involutive of the other one.
To a point $x \in N_o$ corresponds a
quadric $Q(x)$ for which there exists, by construction, only one
family of maximal linear subspaces, and hence this unique family
parametrizes the sections of a
unique sheaves $\E$ corresponding to its unique family of sections and
hence to $x$.
It follows that $\phi^{-1}(N_o)$ is in the ramification locus of
the above covering which coincides with the fixed locus of the
above involution.
\end{proof}

Using the same procedure as in previous proposition we show 
that the same result holds for the base locus of an Enriques net of quadrics.
In this case the discriminant of the net splits as
two smooth cubics intersecting transversally.

\par
\smallskip
{\it Proof of Theorem \ref{thm:reye}.}
Let us consider the case of a K3 $S$ base locus of an Enriques
net of quadrics $N$.
Also in this case,
as in the general one, we have
a 2:1 map $\M_S \to N$
ramified along a plane sextic $N_o \subset N$
relating the moduli space with the
net of quadrics and hence with the K3 surface $S$.
The difference with respect to the general case is that in this case the plane sextic,
which is the discriminant of the net, is
split as the union of two plane cubics $\Delta_A$ and $\Delta_B$.
For,
singular points of $N_o$ correspond
to rank--k quadrics contained in the net with
vertex at least of dimension one.
Since the dimension of the vertex is equal to
$5-k$, $Q$ is a singular point for $N_o$
if $5-k \geq 1$ hence if ${\rm rk}(Q)=k \leq 4$.
Since for a general Enriques net of quadrics
there are no quadrics containing $S$ of rank minor then or equal to $3$ 
we have only
the rank--4 ones to consider.
Note that on $\M_S$ we have two possible involutions to consider;
the first one induced by the covering $\pi$
which interchanges the sheets of the cover, i.e. the covering involution,
and the
other one induced by the involution $i_S$.
These are equivalent and the fixed locus of both is the preimage
of the plane sextic $N_o=\Delta_A \cup \Delta_B$, as in the general case.
When $T$ is of type Reye Congruence we claim the following:\par
{\it Claim: $\M_S$ is a partner of $S$.}
The claim follows by the proof of Theorem 6.1. in \cite{morrison-saito},
since $\M_S$ has the same period of $S$.

\begin{example} \label{ex:(5,16)}
\normalfont{
Let us consider the case of a Kummer surface $S$.
The unicity of the partner is still known by \cite{Muk3}
for reason of the Picard number.
On the other hand, setting $N=NS(S)$, since the map
$O(N) \to O(A_{N})$ is surjective 
(see \cite{morrison-saito})
and since the determinant of $N$
is equal to $2^6$, then by Corollary 3.7 in \cite{HLYO2}
one gets the required unicity result.
}
\end{example}

\

\

\

Address: Dottorato di Ricerca in Matematica,
Dipartimento di Matematica, Universit\`a degli Studi di 
Roma "Tor Vergata", Via della
Ricerca Scientifica, 00133 Roma, Italia

\

email: madonna@mat.uniroma3.it


\begin{thebibliography}{ASDE}

\bibitem{BPV} W.Barth, C.Peters and A.Van de Ven,
Compact complex surfaces, Springer-Verlag, 1984.

\bibitem{bridgeland-maciocia} T.Bridgeland and A.Maciocia,
\emph{Complex surfaces with equivalent derived categories},
Math. Z. 236 (2001), no. 4, 677--697.

\bibitem{caldararu} A.H.Caldararu, \emph{Non--fine moduli spaces
of sheaves on K3 surfaces}, to appear in IMRN, 
preprint math.AG/0108180, 27 august 2001.

\bibitem{cassels} J.W.S.Cassels, Rational quadratic forms,
Academic Press, 1978.

\bibitem{cossec} F.R.Cossec, \emph{Reye Congruences}, Trans.~Amer.~Math.~Soc.
\textbf{280} no.2 (1983), 737-751.

\bibitem{HLYO1} S.Hosono, B.H.Lian, K.Oguiso, and S.T.Yau,
\emph{Autoequivalences of derived category of a K3 surface and monodromy
transformations}, preprint math.AG/0201047, v.2 16 Jan 2002.

\bibitem{HLYO2} S.Hosono, B.H.Lian, K.Oguiso, and S.T.Yau,
\emph{Counting Fourier--Mukai partners and deformations}, preprint math.Ag/0202014,
5 Feb. 2002.

\bibitem{madonna-retta} C.Madonna, \emph{Isogenous K3 surfaces}, 
PhD Thesis, in progress.

\bibitem{miranda-morrison} R.Miranda and D.Morrison,
\emph{The number of embeddings of integral quadratic forms, I}, Proc. Jap. Acad.
\textbf{61}(1985), Ser. A, no.10, 317--320.

\bibitem{morrison-todorov} D.R.Morrison, \emph{On the moduli of Todorov surfaces},
Algebraic Geometry and Commutative Algebra in Honor of Masayoshi Nagata
(H. Hijikata et al., eds.), vol. \textbf{1}, Kinokuniya, Tokyo, 1988, pp. 313-355.

\bibitem{morrison} D.R.Morrison, \emph{Algebraic cycles on products of surfaces},
Proc. Algebraic Geometry Symposium, T${\rm \hat{o}}$hoku University, 1984, pp. 194-210.

\bibitem{morrison-saito} D.R.Morrison and M.H.Saito,
\emph{
Cremona transformations and degrees of period maps for K3 surfaces
with ordinary double points}, Algebraic Geometry, Sendai 1985 (T. Oda, ed.),
Adv. Studies in Pure Math., vol. \textbf{10}, pp. 477-513.

\bibitem{Muk3} S.Mukai, \emph{On the moduli space of bundles
on $K3$ surface, I}, Vector bundles on
algebraic varieties (Bombay, 1984), 341--413,
Tata Inst. Fund. Res. Stud. Math., 11, Tata Inst. Fund. Res.,
Bombay, 1987.

\bibitem{nikulin} V.V.Nikulin, \emph{Integral
symmetric bilinear forms and some of their
applications}, Math. USSR Izv. \textbf{14}(1980),
103--167.

\bibitem{nikulin-zero} V.V.Nikulin,
\emph{On correspondences between K3 surfaces}
Math. USSR Izvestiya, \textbf{30}(1988), no.2, 377--



\bibitem{oguiso} K.Oguiso, \emph{K3 surfaces via almost prime},
math.AG/0110282, to appear on Math.Res.Letter.

\bibitem{Orlov} D.O.Orlov, \emph{Equivalences of derived categories and $K3$
surfaces. Algebraic Geometry,7}, J.Math.Sci.(New York) \textbf{84},
no.5, 1361--1381.

\bibitem{orlov-abelian} D.O.Orlov, \emph{On equivalences
of derived categories of coherent sheaves
on abelian varieties}, alg-geom 9712917, 16 Dec. 1997.

\bibitem{Verra2} A.Verra, \emph{The \'Etale double covering of an
Enriques' surface}, Rend. Sem. Mat. Univers. Politecn. Torino \textbf{41} n.3 (1983),
131--167.

\end{thebibliography}
\end{document}